\documentclass{proc-l}

\usepackage{hyperref, amsrefs,  amssymb,amsbsy,amsfonts,amsthm,latexsym,amsopn,amstext, amsxtra,euscript,amscd}

\newtheorem{thm}{Theorem}
\newtheorem{lem}{Lemma}

\newtheorem{prop}{Proposition}

\newtheorem{exa}{Example}

\newcommand{\set}[1]{{\left\{#1\right\}}}
\newcommand{\K}{{K}}
\def\C{\mathbb C}

\def\F{\mathbb F}

\def\Z{\mathbb Z}
\def\O{\mathcal O}
\def\R{\mathcal R}

\def\a{\alpha}
\def\La{\Lambda}
\def\Th{\Theta}
\def\th{\theta}

\def\om{\omega}

\def\<{\langle}
\def\>{\rangle}

\def\O{\mathcal O}

\def\L{\Lambda}

\numberwithin{equation}{section}

\begin{document}

\title[Codes and theta functions]{Codes over rings of size four, Hermitian lattices, and corresponding theta functions}

\author{T. Shaska}
\address{Department of Mathematics and Statistics,  Oakland University, 368 Science and Engineering Building,
Rochester, MI, 48309.}
\email{shaska@oakland.edu}

\author{G. S. Wijesiri}
\address{Department of Mathematics and Statistics,  Oakland University, 368 Science and Engineering Building,
Rochester, MI, 48309.} \email{gwijesi@oakland.edu}

\date{}




\maketitle
\begin{abstract}
Let $K=Q(\sqrt{-\ell })$ be an imaginary quadratic field with ring of integers $\O_K$, where $\ell$ is a
square free integer such that $\ell\equiv 3 \mod 4$ and $C=[n, k]$ be a linear code defined over
$\O_K/2\O_K$. The level $\ell$ theta function $\Th_{\L_{\ell} (C) } $ of $C$ is defined on  the lattice
$\L_{\ell} (C):= \set { x \in \O_K^n : \rho_\ell (x) \in C}$, where $\rho_{\ell}:\O_K \rightarrow \O_K/2\O_K$
is the natural projection.   In this paper, we prove that:
i) for any $\ell, \ell^\prime$ such that $\ell \leq \ell^\prime$,   $\Th_{\Lambda_\ell}(q)$ and
$\Th_{\Lambda_{\ell^\prime}}(q)$ have the same coefficients up to $q^{\frac {\ell+1}{4}}$,
ii)  for $\ell \geq \frac {2(n+1)(n+2)}{n} -1$,   $\Th_{\L_{\ell}} (C)$ determines the code $C$ uniquely,
iii) for $\ell < \frac {2(n+1)(n+2)}{n} -1$ there is a positive dimensional family of symmetrized weight
enumerator polynomials corresponding to $\Th_{\La_\ell}(C)$.
\end{abstract}


\section{Introduction}
Let $K=Q(\sqrt{-\ell })$ be an imaginary quadratic field with ring of integers $\O_K$, where $\ell$ is a
square free integer such that $\ell\equiv 3 \mod 4$. Then the image $\O_K/2\O_K$ of the projection $
\rho_{\ell}:\O_K \rightarrow \O_K/2\O_K $ is $\F_4$  (resp., $\F_2 \times \F_2$) if $\ell \equiv 3 \mod 8$
(resp.,  $\ell \equiv 7 \mod 8$).

Let $\R$  be a ring isomorphic to $\F_4$  or $\F_2 \times \F_2$ and $C=[n, k]$ be a linear code over $\R$ of
length $n$ and dimension $k$. An admissible level $\ell$ is an $\ell$ such that  $\ell \equiv 3 \mod 8$ if
$\R$ is isomorphic to $\F_4$ or $\ell \equiv 7 \mod 8$ if $\R$ is isomorphic to $\F_2 \times \F_2$. Fix an
admissible $\ell$ and define $\L_{\ell} (C):= \set { x \in \O_K^n : \rho_\ell (x) \in C}.$
Then, the \textbf{level $\ell$ theta function} $\Th_{\L_{\ell} (C) } (\tau) $ of the lattice $\L_{\ell} (C)$
is given as the symmetric weight enumerator $swe_C $ of $C$, evaluated on the  theta functions defined on
cosets of $\O_K/2\O_K$. In this paper we study the following two questions:

i) How do the theta functions $\Th_{\L_{\ell} (C) } (\tau) $ of the same code $C$ differ for different levels
$\ell$?

ii) Can non-equivalent codes give  the same theta functions for all levels $\ell$?

\noindent In an attempt to study the second question Chua in \cite{ch} gives  an example of two
non-equivalent codes that give the same theta function for  level $\ell=7$ but not for higher level thetas.
We will show in this paper how such an example is not a coincidence. Our main results are as follows:

\noindent \textbf{Theorem 1:}   \emph{Let $C$ be a code defined over $R$. For all admissible $\ell,
\ell^\prime$ such that $\ell > \ell^\prime,$ the following holds
\begin{equation*}
\Th_{\Lambda_\ell}(C)=\Th_{\Lambda_{\ell^\prime}}(C)+ {\mathcal O} (q^{\frac {\ell^\prime+1}{4}}).
\end{equation*}}

\noindent \textbf{Theorem 2:} \emph{Let $C$ be a code of size $n$ defined over $\R$ and $\Th_{\La_\ell} (C)$
be its corresponding theta function for level $\ell$. Then the following hold:}

\begin{description}
\item [i)]
\emph{For $\ell < \frac {2(n+1)(n+2)}{n} -1$ there is a $\delta$-dimensional family of symmetrized weight
enumerator polynomials corresponding to $\Th_{\La_\ell}(C)$, where\\ $\delta \geq \frac
{(n+1)(n+2)}{2}-\frac{n(\ell+1)}{4} - 1$.}

\item [ii)]
\emph{For $\ell \geq \frac {2(n+1)(n+2)}{n} -1$ and $n < \frac{\ell+1}{4}$ there is a unique symmetrized
weight enumerator polynomial which corresponds to $\Th_{\La_\ell}(C)$.}
\end{description}

This paper is organized as follows. In the second section, we give a basic introduction of lattices and theta
functions. We define a lattice $\La$ over a number field $K$ in general, the theta series of a lattice, and
one dimensional theta series and its shadow. Then we discuss the  lattices over imaginary quadratic fields
$K=Q(\sqrt{-\ell })$ with ring of integers $\O_K$, where $\ell$ is a square free integer such that
$\ell\equiv 3 \mod 4$. The ring $\O_K / (2\O_K)$ is equivalent to either the field of order 4 or a ring of
order 4 depending on whether $\ell\equiv 3 \mod 8$ or $\ell\equiv 7 \mod 8$. We define bi-dimensional theta
functions for the four cosets of $\O_K / (2\O_K)$.

In the third section, we define codes over $\F_4$ and $\F_2\times \F_2$,  the weight enumerators of a code,
and recall the main result of \cite{ch}. We simplify the expressions for bi-dimensional theta series and
prove Theorem 1.

In the fourth section, we study families of codes corresponding to the same theta function. We call an
\textbf{acceptable theta series} $\Th (q)$ a series for which there exists a code $C$ such that $\Th (q) =
\Th_{\Lambda_\ell}(C) (q)$. For any given $\ell$ and an acceptable theta series $\Th (q) $ we can determine a
family of symmetrized weight enumerators which correspond to $\Th (q)$. For small $\ell$ this is a positive
dimensional family, where the dimension is given by Theorem 2, i). Hence, the example given in \cite{ch} is
no surprise. For large $\ell$ (see Theorem 2, ii)) this is a 0-dimensional family of symmetrized weight
enumerators which correspond to $\Th (q)$. Therefore, the example that Chua provides can not occur for larger
$\ell$.

\section{Introduction to latices and theta functions}
Let $K$ be a number field and $\O_K$ be its ring of integers. A  lattice $\La$ over $K$ is an
$\O_K$-submodule of $K^n$ of full rank. The Hermitian dual is defined by
\begin{equation}
    \L^* = \set{x \in \K^n \; | \; x \cdot \bar{y} \in \O_K, \textit{for all }  y \in \La},
\end{equation}
where $ x \cdot y := \sum_{i=1}^n x_i y_i$. In the case that $\L$ is a free $\O_K$ - module, for every $\O_K$
basis $ \set {v_1 , v_2 , ...., v_n}$ we can associate a Gram matrix G($\L$) given by  $G(\L) =
(v_i.v_j)_{i,j=1}^n $ and the determinant $\det \L :=\det(G) $ defined up to squares of units in $\O_K$. If $
\L = \L^* $ then $ \L$ is Hermitian self-dual (or unimodular) and integral if and only if $\L \subset \L^*$.
An integral lattice has the property $   \L \subset \L^* \subset \frac{1}{det \L} \L$. An integral lattice is
called even if $x\cdot x\equiv 0 \mod 2$ for all $x\in\L$, and otherwise it is odd. An odd unimodular lattice
is called a Type 1 lattice and even unimodular lattice is called a Type 2 lattice.

The theta series of a lattice $\L$ in $K^n$ is given by $ \Th_\L(\tau) = \sum_{z \in \L} e^{\pi i\tau
z\bar{z}}$,  where $ \tau \in H =\set { z \in \C : Im(z)>0}.$   Usually we let $q = e^{\pi i \tau}. $ Then, $
\Th_\L(q) = \sum_{z \in \L} q^{z\bar{z}}$.  The 1-dimensional theta series and its \textbf{shadow} are given
by
\begin{equation}
 \th_3(q):= \sum_{m \in \Z} q^{m^2}, \quad  \th_2(q):= \sum_{m \in {\Z+1/2}} q^{m^2}.
\end{equation}
Let $\ell>0$ be a square free integer and $K=Q(\sqrt{-\ell })$ be the imaginary quadratic field with
discriminant $d_K$. Recall that $d_K=   -\ell$ if $\ell\equiv 3 \mod 4$ and $d_K=-4\ell$ otherwise.

Let $ \O_K$ be ring of integers of $K$. The Hermitian lattice $\L$ over $K$ is an $\O_K $ - submodule of
$K^n$ of full rank. Let $\ell\equiv 3 \mod 4$ and $d$ be a positive number such that $\ell=4d-1$. Then,
$-\ell\equiv 1 \mod 4$. This implies that the ring of integers is $\O_K=\Z[\om_\ell]$, where
$\om_\ell=\frac{-1+\sqrt{-\ell}}{2}$ and $\om_\ell^2 + \om_\ell+d=0$. The principal norm form of $K$ is given
by $ Q_d(x,y) = |x-y\om_\ell|^2 = x^2+xy+dy^2$. Since  $\ell\equiv 3 \mod 4$, we can consider two cases:

(1) If $\ell \equiv 3 \mod 8$ then $-\ell \equiv 5 \mod 8$. Thus, the prime ideal $\<2\> \subset \Z$ lifts to
a prime $2\O_K \subset \O_K$.  Since the ring of integers $\O_K$ is a Dedekind domain, $2\O_K$ is a maximal
ideal. Therefore $\O_K / (2\O_K) \simeq \F_4$.

(2) If $\ell \equiv 7 \mod 8$ then $-\ell \equiv 1 \mod 8$. Then the prime ideal $\<2\> \in \Z$ splits in
$K$. Therefore $2\O_K$ splits in $\O_K$. Hence, $\O_K / (2\O_K) \simeq \F_2 \times \F_2$. In either case, a
complete set of coset representatives is $\set{0,1,\om_\ell,1+\om_\ell}$.

Let the following be the bi-dimensional theta series for the four cosets:
\begin{equation}
\begin{split}
A_d(q) & := \Th_{2\O_K}(\tau) =\sum_{m,n \in \Z} q^{4Q_d(m,n)} \\
C_d(q) & := \Th_{1+2\O_K}(\tau) =\sum_{m,n \in \Z} q^{4Q_d(m+\frac{1}{2},n)} \\
G_d(q) & := \Th_{\om_\ell+2\O_K}(\tau) =\sum_{m,n \in \Z} q^{4Q_d(m,n+\frac{1}{2})}\\
H_d(q) & := \Th_{1+\om_\ell+2\O_K}(\tau) =\sum_{m,n \in \Z} q^{4Q_d(m+\frac{1}{2},n+\frac{1}{2})}
\end{split}
\end{equation}
Then we have the following lemma.
\begin{lem}
Bi-dimensional theta series can be further expressed in terms of the standard one dimensional theta series
and its shadow.
\begin{equation}\label{eq2}
\begin{split}
A_d(q)&= \th_3(q^{4}) \th_3(q^{4\ell}) + \th_2(q^{4}) \th_2(q^{4\ell}) \\
C_d(q)&= \th_2(q^{4}) \th_3(q^{4\ell}) + \th_3(q^{4}) \th_2(q^{4\ell}) \\
G_d(q)&=H_d(q)=\frac{ \th_2(q) \th_2(q^\ell)}{2}.
\end{split}
\end{equation}
Moreover,
\begin{equation}\label{eq3}
 2G_d(q)=A_d(q^{1/4})-A_d(q)-C_d(q).
\end{equation}
\end{lem}
\proof  See \cite{sl1} for details. \endproof

\section{Codes over $\F_4$ and $\F_2\times \F_2$}
Let $\F_4 =\set{0,1,\om,\om^2}$, where $\om^2+\om+1=0$, be the finite field of four elements. The conjugation
is given by $\bar{x}=x^2$, $x \in \F_4$. In particular $\bar{\om}=\om^2=\om+1$. Let $R_4=\F_2+\om\F_2$ with
the new equation for $\om$ is  $\om^2+\om=0$. Notice that $R_4$ has two maximal ideals namely $\<\om\>$ and
$\< \om+1 \>$. Furthermore, one can show that $R_4/\<\om\>$ and $R_4/\<\om+1\>$ are both isomorphic to
$\F_2$. The Chinese remainder theorem tells us that $R_4=\<\om\>\oplus\<\om+1\>$. Therefore, $R_4\simeq \F_2
\times \F_2$. The conjugate of $\om$ is $\om+1$. Let $\R$ be the field $\F_4$ if $ \ell \equiv  3\mod 8$ or
the ring $R_4\simeq \F_2 \times \F_2$ when $ \ell \equiv 7\mod 8 $.  A linear code $C$ of length $n$ over
$\R$ is an $\R$-submodule of $\R^n$. The dual is defined as $C^\bot=\set{u\in \R: u\cdot \bar{v}=0 \textit{
for all } v \in C}$. If $C=C^\bot$ then $C$ is self-dual.

We define $\L_{\ell}(C):= \set { x \in \O_K^n : \rho_{\ell}(x) \in C} $ where $\rho_{\ell}:\O_K \rightarrow
\O_K/2\O_k \rightarrow \R$ . In other words, $\L_{\ell}(C)$ consists of all vectors in $\O_K^n$ which when
taken mod $2\O_K$ componentwise are in $\rho_{\ell}^{-1}(C)$. The following is immediate.
\begin{lem}
\begin{enumerate}
\item  $\L_{\ell}(C)$ is an $\O_K$-lattice.
\item $\L_{\ell}(C^\bot)=2\L_{\ell}(C)^*.$
\item  $C$ is self dual if and only if $\frac{\L_\ell(C)}{\sqrt{2}} $ is self dual.
\end{enumerate}
\end{lem}
Let $u=(u_1,u_2,\cdots,u_n) \in \R^n$ be a codeword and $\a \in \R$. Then the counting function $n_\a(u)$ is
defined as the number of elements in the set  $\set{j:u_j=\a}$. For a code $C$ we define the complete weight
enumerator ($cwe$), symmetrized weight enumerator ($swe$) and Hamming weight enumerator ($W$) to be
\begin{equation}
\begin{split}
cwe_C(X,Y,Z,W)&:=\sum_{u \in C}X^{n_0(u)}Y^{n_1(u)}Z^{n_\om(u)}W^{n_{1+\om}(u)}\\
swe_C(X,Y,Z)&:=cwe_C(X,Y,Z,Z)\\
W_C(X,Y)&:=swe_C(X,Y,Y).\\
\end{split}
\end{equation}
Then we have the  following.
\begin{prop}
Let $\ell\equiv 3\mod 4$, $C$ be a linear code over $\R$, and $\frac{\L_\ell(C)}{\sqrt{2}} $ be a Hermitian
lattice constructed via the  construction A. Then
\begin{equation}\label{eq4}
\theta_{\L_\ell(C)}(\tau)=swe_C(A_d(q),C_d(q),G_d(q))
\end{equation}
where $A_d(q)$, $C_d(q)$, and $G_d(q)$ are given as in  Eq.~\eqref{eq2}.
\end{prop}

For a proof of the above statement the reader can see \cite{ch}. From the definition of one dimensional theta
series we have
 \[
 \theta_2(q)=2q^{1/4}\sum_{i\in S}q^i, \quad \theta_2(q^4)=2q\sum_{i:odd}q^{i^2-1}, \quad \theta_3(q^4)=1+2q^4\sum_{i\in\Z^+}q^{4(i^2-1)},
\]
where $S=\set{ \frac{j^2-1}{4}: j\equiv 1\mod 2}$. From Eq.~\eqref{eq2} we can write
\[G_d(q) =\frac{\th_2(q)\th_2(q^{\ell})}{2} = q^{\frac{(\ell+1)}{4}}\alpha_1,\] where $\alpha_1=\sum_{i \in S}q^i\sum_{j \in S}q^{\ell j}$.
Then,
\[\begin{split}
 A_d(q) & = \th_3(q^{4}) \th_3(q^{4\ell}) + \th_2(q^{4}) \th_2(q^{4\ell})\\
              & = (1+2q^4\sum_{i \in \Z^+}q^{4(i^2-1)}) (1+2q^{4\ell}\sum_{j \in \Z^+}q^{4\ell(j^2-1)})\\
              & + 4q^{\ell+1}\sum_{i:odd}q^{i^2-1} \sum_{j:odd}q^{(j^2-1)\ell}\\
              & = \alpha_2+q^{\ell+1}\alpha_3+q^{4\ell}\alpha_4,
\end{split}
\]
where $\alpha_2, \alpha_3$ and $\alpha_4$ have the following forms
\[
\begin{split}
 \alpha_2 &=1+2q^4\sum_{i \in \Z^+}q^{4(i^2-1)}\\
 \alpha_3 &=4\sum_{i:odd}q^{i^2-1} \sum_{j:odd}q^{(j^2-1)\ell}\\
 \alpha_4 &=2\sum_{j \in\Z^+}q^{4\ell(i^2-1)}(1+2q^4\sum_{i \in\Z^+}q^{4(i^2-1)}).
\end{split} \]
\noindent Furthermore,
\[\begin{split}
 C_d(q) & = \th_2(q^{4}) \th_3(q^{4\ell}) + \th_3(q^{4}) \th_2(q^{4\ell})\\
              & = 2q\sum_{i:odd}q^{i^2-1} (1+2q^{4\ell}\sum_{i \in \Z^+}q^{4\ell(i^2-1)})\\
              & + (1+2q^4\sum_{i \in \Z^+}q^{4(i^2-1)})(2q^{\ell} \sum_{i:odd}q^{(i^2-1)\ell}\\
              & = \alpha_5+q^{\ell}\alpha_6+q^{4\ell+1}\alpha_7,
\end{split}
\]
where $\alpha_5, \alpha_6$ and $\alpha_7$ have the following forms
\[\begin{split}
  \alpha_5 &=2\sum_{i:odd}q^{i^2-1}\\
  \alpha_6 &=2\sum_{j:odd}q^{(j^2-1)\ell}(1+2q^4\sum_{i \in \Z^+}q^{4(i^2-1)})\\
  \alpha_7 &=4\sum_{i:odd}q^{i^2-1}\sum_{j \in \Z^+}q^{4\ell(j^2-1)}.
\end{split}\]
The next result shows  that for large enough admissible $\ell$ and $\ell^\prime$  the theta functions
$\Th_{\Lambda_\ell}(C)$ and $\Th_{\Lambda_{\ell^\prime}}(C)$ are virtually the same.
\begin{thm}
Let $C$ be a code defined over $R$. For all admissible $\ell, \ell^\prime$ such that $\ell > \ell^\prime,$
the following holds
\begin{equation}
\Th_{\Lambda_\ell}(C)=\Th_{\Lambda_{\ell^\prime}}(C)+ {\mathcal O} (q^{\frac {\ell^\prime+1}{4}}).
\end{equation}
\end{thm}
\proof Let \[swe_C(X, Y, Z) =\sum_{i+ j+ k = n} a_{i,j,k} \cdot X^i Y^j Z^k\] be a degree $n$ polynomial.
Write this as a polynomial in $Z$. Then \[swe_C(Z)=\sum_{k=0}^n L_kZ^k = L_0 + Z(\sum_{k=1}^n L_kZ^{k-1}).\]
Terms in $L_0$ are of the form of $a_{i,j} X^i Y^j $, where $i+j=n$. From the above we have
\[
\begin{split}
A_d(q)^i\cdot C_d(q)^j &=(\alpha_2+q^{\ell+1}\alpha_3+q^{4\ell}\alpha_4)^i\cdot(\alpha_5+q^{\ell}\alpha_6+q^{4\ell+1}\alpha_7)^j\\
                      &=(\textit{terms independent from} ~ \ell)+q^\ell(\cdots)
\end{split}
\]
Also we have seen that $G_d(q)=q^{(\ell+1)/4}\alpha_1$. This gives
\[
\begin{split}
\Th_{\Lambda_\ell}(C)&=swe_C(A_d(q),C_d(q),G_d(q))\\&=(\textit{terms independent from}~ \ell)+{\mathcal O}
(q^{\frac {\ell+1}{4}}) \end{split}. \] Then the result follows.
\endproof

\begin{exa}
Let $C$ be a code defined over $R_4$ which has symmetrized weight enumerator
\[swe_{C}(X,Y,Z) = X^3+X^2Z+XY^2+2XZ^2+Y^2Z+2Z^3. \]
Then we have the following:
\begin{equation}
\begin{split}
\Th_{\Lambda_{63}}(C) &= 1 + 6q^4 + 12q^8 + 8q^{12} + 12q^{16} + 6q^{18} + 48q^{20} + 30q^{22} + \cdots\\
\Th_{\Lambda_{79}}(C) &= 1 + 6q^4 + 12q^8 + 8q^{12} + 6q^{16} + 30q^{20} + 6q^{22} + 48q^{24} + \cdots\\
\Th_{\Lambda_{79}}(C) &= \Th_{\Lambda_{63}}(C)+ {\mathcal O} (q^{16}).
\end{split}
\end{equation}
\end{exa}

\section{A family of codes corresponding to the same theta function}
If we are given the code over $\R$ and its symmetrized weight enumerator polynomial, then by Eq.~\eqref{eq4}
we can find the theta function of the lattice constructed from the code by using the construction \emph{A}.
Now, we would like to give a way to construct families of codes corresponding to the same theta function.

Let $\Th(q) = \sum_{i=0}^\infty \lambda_i q^i$ be an acceptable theta series for level $\ell$ and
\[f(x,y,z)=\sum_{i+j+k=n}    c_{i,j,k}x^iy^jz^k\]
be a degree $n$ generic ternary homogeneous polynomial. We want to find out how many polynomials $f(x,y,z)$
correspond to $\Th(q)$ for a fixed $\ell$.

We have the following lemma.
\begin{lem}
Let $C$ be a code of size $n$ defined over $\R$ and $\Th(q)$ be its theta function for level $\ell$. Then,
$\Th(q)$ is uniquely determined by its first $\frac{n(\ell+1)}{4}$ coefficients.
\end{lem}

\proof Let $C$ be a  code over $\R$,    $\Th(q) = \sum_{i=0}^\infty \lambda_i q^i$ be its theta series,
$s=\frac{n(\ell+1)}{4}$ and \[f(x,y,z)=\sum_{i+j+k=n} c_{i,j,k}x^iy^jz^k\] be a degree $n$ generic ternary
homogeneous polynomial. Find $A_d(q), C_d(q),$$ G_d(q)$ for the given $\ell$ and substitute in $f(x,y,z)$.
Hence $f(x,y,z)$ is now written as a series in $q$. Recall that a generic degree $n$ ternary polynomial has
$r=\frac{(n+1)(n+2)}{2}$ coefficients. So,  the corresponding coefficients of the two sides of the equation
are equal:
\[ f(A_d(q), C_d(q), G_d(q) ) = \sum_{i=0}^\infty \lambda_i q^i.\]
Consider the term
\[c_{i,j,k}(\alpha_2+q^{\ell+1}\alpha_3+q^{4\ell}\alpha_4)^i(\alpha_5+q^{\ell}\alpha_6
+q^{4\ell+1}\alpha_7)^j(q^{\frac{(\ell+1)}{4}}\alpha_1)^k.\]
Then $c_{i,j,k}$ appears first as a coefficient of $q^{j+\frac{k(\ell+1)}{4}}$. For all such $j,k$ we have
$j+\frac{k(\ell+1)}{4} \leq \frac{n(\ell+1)}{4}$. Consider the equations where $c_{i,j,k}$ appears first.
This is a system of equations with $\leq \frac{(n+1)(n+2)}{2}$ equations.
Let us denote this system of equations as $\Xi$.
Solve this system for $c_{i,j,k}$. Hence,  $c_{i,j,k}$ is a function of
$\l_0,\cdots,\l_s$.  For each $\mu>s, \l_\mu$ is a function of $c_{i,j,k}$ for $i,j,k=0,\cdots,n$, and
therefore a rational function on $\l_0,\cdots,\l_s$. This completes the proof.
\endproof

Next we have the following theorem:
\begin{thm}
Let $C$ be a code of size $n$ defined over $\R$ and $\Th_{\La_\ell} (C)$ be its corresponding theta function
for level $\ell$. Then the following hold:
\begin{description}
\item [i)] For $\ell < \frac {2(n+1)(n+2)}{n} -1$ there is a $\delta$-dimensional family of symmetrized weight
enumerator polynomials corresponding to $\Th_{\La_\ell}(C)$, where\\ $\delta \geq \frac
{(n+1)(n+2)}{2}-\frac{n(\ell+1)}{4} - 1$

\item [ii)]For $\ell \geq \frac {2(n+1)(n+2)}{n} -1$ and $n < \frac{\ell+1}{4}$ there is a unique symmetrized weight
enumerator polynomial which corresponds to $\Th_{\La_\ell}(C)$.
\end{description}
\end{thm}

\proof We want to find out  how many polynomials $f(x,y,z)$ correspond to $\Th_{\La_\ell}(C)$ for a fixed
$\ell$.  $\Th_{\La_\ell}(C) $ and $f(x,y,z)$ are defined as above. Consider the system of equations $\Xi$.

If $\frac{n(\ell+1)}{4} < r$ then our system has more variables than equations. Since the system is linear,
the solution space is a family of positive dimension.

If $\frac{n(\ell+1)}{4} \geq r$ then for each equation in $\Xi$ (see the proof of the previous Lemma) we have
only one $c_{i,j,k}$ appearing for the first time. Otherwise suppose $c_{i,j,k}$ and
$c_{i^\prime,j^\prime,k^\prime}$ appear for the first time in an equation of $\Xi$. Then $ j +
\frac{k(\ell+1)}{4} = j^\prime + \frac{k^\prime(\ell+1)}{4}.$ This implies
\begin{equation}\label{eq8}
4(j-j^\prime)=(k^\prime-k)(\ell+1).\end{equation} Without loss of generality, assume $k^\prime \geq k.$ We
can consider three cases.

case 1: If $k^\prime-k \geq 2$, then from Eq.~\eqref{eq8} we have $4n(j-j^\prime)=n(k^\prime-k)(\ell+1) \geq
4r(k^\prime-k).$ Then we have $n(j-j^\prime) \geq (n+1)(n+2).$ Since $n \geq (j-j^\prime),$ we have a
contradiction.

case 2: If $k^\prime-k =1$, then by Eq.~\eqref{eq8} $j-j^\prime = \frac {\ell+1}{4}.$ Since $j-j^\prime \leq
n$ and $\frac{\ell+1}{4} > n,$ we get a contradiction.

case 3: $k^\prime-k = 0$, then by Eq.~\eqref{eq8} we have $j=j^\prime.$ Hence $i=i^\prime.$

Notice that $c_{n,0,0}$ is uniquely determined by the equation corresponding to the equation of coefficient
of $q^0$. Solve the system $\Xi$ in the order of the equation that corresponds to the power of $q$. We have a
unique solution for $c_{i,j,k}$.
\endproof

\subsection{Families of codes of length 3}
In this section we discuss the codes of length 3 for different levels $\ell$. Our main goal is to investigate
the example provided in \cite{ch} and provide some computational evidence for the above two cases. We assume
that the symmetrized weight enumerator polynomial is a generic homogenous polynomial of degree three.

Let $P(x, y, z)$ be a generic ternary cubic homogeneous polynomial given as below
\begin{equation}
\begin{split}
P(x,y,z)& =c_1x^3+c_2y^3+c_3z^3+c_4x^2y+c_5x^2z+c_6y^2x+c_7y^2z \\
        & +c_8z^2x+c_9z^2y+c_{10}xyz.
\end{split}
\end{equation}

Assume that there is a code $C$, of length 3, defined over $\R$ such that $swe_C (x, y, z)=P(x, y, z)$. First
we have to fix the level $\ell$. When we fix the level, we can find $A_d(q),C_d(q),G_d(q)$. By equating both
sides of \[ p(A_d(q), C_d(q), G_d(q) ) = \sum_{i=0}^\infty \lambda_i q^i,\] we can get a system of equations.
When $\ell=7,$  we are in the first case of the previous theorem. The system of equations  is given by the
following.
\begin{equation}\label{eq9}
\left\{
\begin{aligned}
 & c_1 - \lambda_0 = 0\\
 &2c_4- \lambda_1 = 0\\
 &4c_6+2c_5-\lambda_2 = 0 \\
 &8c_2+4c_{10}- \lambda_3 = 0 \\
\end{aligned}
\right. \qquad \left\{
\begin{aligned}
 &6c_1+4c_8+2c_5+8c_7-\lambda_4 = 0 \\
 &8c_4+8c_9+4c_{10}-\lambda_5 = 0 \\
 &8c_5+8c_3+8c_7+8c_8+8c_6-\lambda_6 = 0.\\
\end{aligned}
\right.
\end{equation}
The solution for the above system is given by $c_1  = \lambda_0, c_4 = \frac{1}{2} \lambda_1$, and
\begin{equation}\label{eq5}
\begin{split}
c_2 &= \frac{1}{2} \lambda_1 + \frac{1}{8} \lambda_3  -\frac{1}{8} \lambda_5 + c_9, \quad c_3 = \frac{3}{2}
\lambda_0 - \frac{1}{4} \lambda_2  - \frac{1}{4} \lambda_4 + \frac{1}{8} \lambda_6 + c_7,
\quad \\
c_5 &= -3 \lambda_0 +\frac{1}{2} \lambda_4  - 4c_7 - 2c_8, \quad c_6 = \frac{3}{2} \lambda_0 + \frac{1}{4}
\lambda_2 - \frac{1}{4} \lambda_4 + 2c_7 + c_8, \\
  c_{10} & = - \lambda_1 + \frac{1}{4} \lambda_5 - 2c_9
\end{split}
\end{equation}
where $c_7,c_8,c_9$ are free variables. By giving different triples $(c_7, c_8, c_9)$, we can construct
different polynomials $P(x,y,z)$ for the same $\sum_{i=0}^\infty \lambda_i q^i.$

Consider the following theta function. From \cite{ch} there are two non isomorphic codes that give this theta
function for level $\ell=7$:
\begin{equation}
\Th_{\sqrt{2}{\O_K}_7^3} = 1+6q^2+24q^4+56q^6+114q^8+168q^{10}+280q^{12}+294q^{14}+\cdots
\end{equation}
For this particular theta function, we can rewrite the solution (Eq. \eqref{eq5}) as follows:
$c_1=1,c_2=c_9,c_3=1+c_7,c_4=0, c_5=9-4c_7-2c_8,c_6=-3-2c_7+c_8,c_{10}=-2c_9.$

For the triple $(1,2,0)$ (resp., $(0,3,0)$) we get the symmetrized weight enumerator polynomial for the code
$C_{3,2}$ (resp. $C_{3,3}$). That is $ swe_{C_{3,2}}(X,Y,Z)=X^3+X^2Z+XY^2+2XZ^2+Y^2Z+2Z^3$ (resp.,  $
swe_{C_{3,3}}(X,Y,Z)=X^3+3X^2Z+3XZ^2+Z^3$),
where $C_{3,2}$ and $C_{3,3}$ are given by:
\begin{equation}
\begin{split}
C_{3,2} & = \om<[0,1,1]> + (\om+1)<[0,1,1]>^\bot \\
C_{3,3} & = \om<[0,0,1]> + (\om+1)<[0,0,1]>^\bot.
\end{split}
\end{equation}
When $\ell=15,$ we are in the second case of the above theorem. The system of equations is as follows:
\begin{equation}\label{eq10}
\left\{
\begin{aligned}
 &c_1 - \lambda_0 = 0\\
 &2c_4 - \lambda_1 = 0\\
 &4c_6 - \lambda_2 = 0 \\
 &8c_2 - \lambda_3 = 0 \\
 &6c_1 + 2c_5 - \lambda_4 = 0 \\
\end{aligned}
\right. \qquad \left\{
\begin{aligned}
 &8c_4 + 4c_{10} - \lambda_5 = 0 \\
 &2c_5 + 8c_7 + 8c_6 - \lambda_6 = 0\\
 &4 c_8 + 8c_7 + 12c_1 + 8c_5 - \lambda_8 = 0\\
 &10c_4 + 8c_9 +8c_{10} - \lambda_9 = 0\\
 &8c_7 + 8c_5 + 12c_8 + 8c_3 + 8c_1 - \lambda_{12} = 0. \\
\end{aligned}
\right.
\end{equation}
Each $c_i$ appears first in exactly one equation. For example consider the seventh equation. $c_7$ is the
only variable that appears first in the seventh equation. Solve the system in given order. The solution for
the above system is given by; $c_1 = \lambda_0$, $c_2 = \frac{1}{8} \lambda_3$, $c_4 = \frac{1}{2}
\lambda_1$, $c_6 = \frac{1}{4} \lambda_2$, and
\begin{equation}
\begin{split}
\begin{aligned}
& c_3 = -\lambda_0 - \frac{1}{2} \lambda_2 + \frac{3}{4} \lambda_4 + \frac{1}{4} \lambda_6 - \frac{3}{8}
\lambda_8 +  \frac{1}{8} \lambda_{12},\\
& c_7 = \frac{3}{4} \lambda_0 - \frac{1}{4} \lambda_2  - \frac{1}{8} \lambda_4 + \frac{1}{8} \lambda_6,\\
& c_8 = \frac{3}{2} \lambda_0 + \frac{1}{2} \lambda_2  - \frac{3}{4} \lambda_4 - \frac{1}{4} \lambda_6 +
\frac{1}{4} \lambda_8,\\
\end{aligned}
\quad
\begin{aligned}
&c_5 = -3 \lambda_0 + \frac{1}{2} \lambda_4 \\
& c_9 = \frac{3}{8} \lambda_1 - \frac{1}{4} \lambda_5 + \frac{1}{8} \lambda_9\\
& c_{10} = - \lambda_1 + \frac{1}{4} \lambda_5\\
\end{aligned}
\end{split}
\end{equation}
We have a unique solution. This implies that two non equivalent codes cannot give the same theta function for
$\ell=15$ and $n=3$.

\section{Concluding remarks}
The main goal of this paper was to find out  how theta functions determine the codes over a ring of size 4.
First we have shown  how the theta functions of the same code $C$ differ for different levels $\ell$. The
first $\frac{\ell+1}{4}$ terms of the theta functions for levels $\ell$ and $\ell^\prime$ are the same, where
$\ell^\prime \geq \ell$.

In \cite{ch}, two non-isomorphic codes that give the same theta function for level $\ell=7$ but not under
higher level constructions are given. We justified the reason why we don't have a similar situation for
higher level constructions. In this note we have addressed a method that we can use for finding a family of
polynomials that correspond to a given acceptable theta series for some fixed level $\ell$. We have studied
two cases depending upon $\ell$ that give either a positive dimensional family of polynomials or a unique
polynomial.

\end{document}